\documentclass[final,leqno,letterpaper]{amsart}

\usepackage{amsmath}
\usepackage{graphicx,amssymb}

\usepackage{cases}
\usepackage{float}
\usepackage{listings}
\usepackage{times}
\usepackage{tikz}
\usepackage{epic}
\usepackage{geometry}
\usepackage{bm}
\usepackage{stmaryrd}
\usepackage{indentfirst}
\usepackage{xcolor}
\usepackage{epstopdf}
\usepackage{enumerate}
\usepackage{enumitem}
\usepackage{hyperref}

\usepackage[foot]{amsaddr}

\makeatletter
\renewcommand{\email}[2][]{%
  \ifx\emails\@empty\relax\else{\g@addto@macro\emails{,\space}}\fi%
  \@ifnotempty{#1}{\g@addto@macro\emails{\textrm{(#1)}\space}}%
  \g@addto@macro\emails{#2}%
}

\hypersetup{%
    pdftitle={A deep learning based nonlinear upscaling method for transport equations},
    pdfauthor={Tak Shing Au Yeung, Eric T. Chung, Simon See},
    pdfkeywords={a keyword, another keyword, and another one}
    }

\title{A deep learning based nonlinear upscaling method for transport equations
}

\author{Tak Shing Au Yeung$^1$}
\address{$^1$ Department of Mathematics,The Chinese University of Hong Kong, Hong Kong Special Administrative Region.}
\author{Eric T. Chung$^2$}
\address{$^2$ Department of Mathematics,The Chinese University of Hong Kong, Hong Kong Special Administrative Region.}
\author{Simon See$^3$}
\address{$^3$NVIDIA - NVAITC, Santa Clara, United States.}

\begin{document}

\begin{abstract}
We will develop a nonlinear upscaling method for nonlinear transport equation. The proposed scheme gives a coarse scale equation
for the cell average of the solution. In order to compute the parameters in the coarse scale equation, a local downscaling operator is constructed.
This downscaling operation recovers fine scale properties using cell averages. This is achieved by solving the equation on an oversampling region with
the given cell average as constraint. 
Due to the nonlinearity, one needs to compute these downscaling operations on the fly and cannot pre-compute these quantities. 
In order to give an efficient downscaling operation, we apply a deep learning approach. We will use
a deep neural network to approximate the downscaling operation. Our numerical results show that the proposed scheme can achieve a good accuracy
and efficiency.  
\end{abstract}

\maketitle

\section{Introduction}

Multiscale methods or numerical upscaling for linear problems have been widely developed
and are very useful for multiscale problems. Nevertheless, many realistic problems are of nonlinear and multiscale nature. 
For example, the dynamics of
multi-phase flow and transport in heterogeneous media
varies over multiple space and time scales. 
In order to accurately capture the coarse scale dynamics, some types of nonlinear upscaling are necessary. 
There are in literature many linear and nonlinear upscaling techniques, and some of these include
\cite{ab05, egw10,  arbogast02, GMsFEM13, AdaptiveGMsFEM, brown2014multiscale, ElasticGMsFEM, ee03, abdul_yun, ohl12, fish2004space, fish2008mathematical, oz07, matache2002two, apwy07, henning2009heterogeneous, OnlineStokes, chung2017DGstokes,WaveGMsFEM, pwy02, Arbogast_PWY_07, MsDG, fish1997computational,oskay2007eigendeformation,yuan2009multiple,chung2010reduced}.
 Nonlinear
upscaling methods, such as pseudo-relative permeability approach
\cite{cdgw03,KB,BT},
computes nonlinear relative permeability functions based on single cell
two-phase flow computations.
It is known that these nonlinear approaches lack
robustness and they are processes dependent \cite{ed01,ed03}. To overcome
these difficulties, one needs to find better nonlinear upscaling techniques,
and it is the goal of this paper to do this.

There are in literature some nonlinear upscaling approaches. 
One example is nonlinear homogenization
\cite{pankov97, ep03d}, for which local nonlinear problems are solved on each coarse grid block and are used in the 
construction of global coarse grid formulation. 
Nonlinear homogenization is applied to many problems including the p-Laplace equations, pseudo-elliptic equations and
parabolic equations  \cite{ep03a,ep03c,ep03d}.
These methods typically require the assumption of scale separation, and may give limited accuracy for applications that do not
admit such assumption. 
Another nonlinear upscaling approach, that is popular in computational mechanics, is the computational continua framework \cite{fish2010computational}.
These methods, for example \cite{fish2010computational,fafalis2018computational}, 
use nonlocal quadrature to couple the coarse scale system stated on a unions of some disjoint computational unit cells with the aim of solving problems with non-scale-separation heterogeneous media.

The nonlinear upscaling method developed in this paper is motivated by the general framework, called Nonlinear Nonlocal Multi-continua Upscaling (NLNLMC), 
of nonlinear upscaling presented in \cite{chung2018NLNLMC}. 
This framework originates from the Constraint Energy Minimizing Generalized Multiscale Finite Element
Method (CEM-GMsFEM) \cite{chung2018constraint,chung2018fast,chung2018constraintmixed}.
The CEM-GMsFEM gives multiscale basis functions that are localizable even for the case of highly heterogeneous and high contrast media. The idea there is
the use of local spectral problems and an energy minimization principle. The basis functions are constructed by solving problems on some oversampling regions. 
In addition, rigorous convergence analysis is presented and shows that the convergence
depends only on the coarse grid size and is independent of the media. 
However, the degrees of freedoms in CEM-GMsFEM do not have physical meaning. 
In upscaling, one desires unknown variables that have physical meanings such as cell averages. 
With this goal in mind, the Nonlocal Multi-continua Upscaling (NLMC) is introduced \cite{chung2018nlmc}. 
The key idea follows CEM-GMsFEM except that the multiscale basis functions are modified so that the degrees of freedom represent
the average of the solutions. Therefore, the resulting coarse scale equation gives a relation between the cell averages of the solutions on the coarse grid,
and a convergence analysis is given in \cite{zhao2020analysis}. 
We remark that both CEM-GMsFEM and NLMC methods are suitable for linear multiscale problems. 

The framework of NLNLMC extends the concept of NLMC method with the aim of finding nonlinear upscaling for nonlinear problems \cite{chung2018NLNLMC}. 
In this case, one needs to avoid the concept of basis functions, as they are only applicable to linear problems. 
In general, the NLNLMC framework has three important methodological ingredients. 
First, we identify macroscopic variables for each coarse grid block, similar to multicontinua variables  \cite{barenblatt1960basic,lee2001hierarchical, warren1963behavior, panasenko2018multicontinuum}.
Secondly, we will construct local downscaling functions. 
In particular, given the macroscopic variables, we will solve local problems on some oversampling regions to recover fine scale properties. 
It is important that these local problems are solved on oversampling regions, and this allows connectivity of neighbouring macroscopic variables. 
Finally, a global coarse scale equation is obtained by combining all local downscaling functions and using a suitable coarse scale solver. 
The resulting scheme is a coarse scale equation that relates all macroscopic variables. Moreover, the connectivity is nonlocal as the oversampling regions
can be of several coarse grid layers wide.

The goal of this paper is to construct a nonlinear upscaling approach for nonlinear transport equations. We will apply the general concept of NLNLMC
together with deep learning techniques. 
For the macroscopic variable, we will use cell average on coarse element, 
and for the global coarse scale solver, we will use an upwind finite volume scheme. 
The main component of the proposed scheme is the local downscaling functions, and these will give the parameters required in the final coarse grid equation. 
Given a set of macroscopic values, we will solve a local problem on an oversampling region to construct a fine scale downscaling function, whose mean values 
on coarse elements match the given macroscopic values. In general, this is an expensive task as
these local problems are solved on-the-fly when the solution averages are given, and one cannot easily pre-compute these problems. 
This fact motivates the use of deep learning. 
There are in literature some works on using deep neural networks to learn macroscopic parameters in 
coarse scale or reduced order models, see for example \cite{wang2020deep,vasilyeva2020learning,wang2020reduced,zhang2020learning,cheung2020deep}. 
The main idea is to consider the macroscopic variables as input and the downscaling functions or their average values as output. Then suitable deep neural networks
are trained and are used to approximate this expensive procedure. The resulting approach allows the use of deep neural network
to learn the parameters required in the coarse scale equations. 
This can give a significant improvement in the computational times, as we will see
in our numerical simulations. We remark that using deep neural network for reduced models allows a robust learning process
as there are fewer parameters to be learned.

The paper is organized as follows. 
In Section~\ref{sec:prelim}, we will present the problem formulations and some basic notations. 
In Section~\ref{sec:method}, the key elements of the proposed method will be discussed in detail. This includes the construction of the method
and the deep neural network, as well as some implementation details.
Computational results will be presented in Section~\ref{sec:num} to validate our scheme. We will
show the performance of our method, and compare our method with a standard scheme without using nonlinear upscaling. We will also
compare the accuracy and efficiency with and without the use of deep neural networks. 
Finally, a conclusion is given in Section~\ref{sec:conclusion}. 

\section{Preliminaries}\label{sec:prelim}
\subsection{Basic setup}
\label{S:1}
We let $\Omega\subset \mathbb{R}^2$ be a computational domain in two space dimensions and let $T>0$ be a fixed time. 
Our goal is to design a nonlinear upscaling method for the 
following transport equation:
\begin{align*}
\frac{\partial S}{\partial t}+\text{div}\left(v  \lambda(S)\right)&=f, \qquad\quad~ \text{ in } (0,T)\times \Omega, \\
S&=g, \qquad \quad~\text{ on } (0,T) \times \Gamma, \\
S(0,x)&=S_0(x), \qquad \text{ in } \Omega,
\end{align*}
where $\Gamma = \{ x\in \partial \Omega \; : \; v\cdot n <0\}$ is the inflow boundary of $\Omega$, $n$ is the outward unit normal vector of $\partial \Omega$ 
and $f,g,S_0$ are given functions. Motivated by applications, 
we assume that the velocity $v$ is divergence free, that is, $\text{div}(v)=0$. 
Moreover, the function $\lambda : \mathbb{R} \rightarrow \mathbb{R}$ is a nonlinear function.
This problem is motivated by two phase flow and transport problems, in which the velocity $v$ is given by the Darcy's law \cite{fu2019local}. In general,
the transport equation and the Darcy's law are coupled. In this paper, we focus only solving the nonlinear transport equations. 
The development of upscaling methods for the coupled problem will be considered in a forthcoming paper. 

Next, we introduce the notions of coarse grids. Let $\mathcal{T}_H$ be a conforming partition of $\Omega$ into  finite elements. Here, $H$ is the coarse mesh size and this partition is called coarse grid. We let $N_K$ be the number of elements in the coarse mesh. Then we assume that each coarse element is partitioned into a connected union of  fine-grid blocks and this partition is called $\mathcal{T}_h$. Note that $\mathcal{T}_h$ is a conforming refinement of the coarse grid $\mathcal{T}_H$ with the mesh size $h$. That is, for any $\tau\in \mathcal{T}_h$, there exist $K\in \mathcal{T}_H$ such that $\tau\subset K$. 
See Figure~\ref{fig:testfig} for an illustration.
In our method, we will develop a nonlinear upscaling technique that gives
the coarse grid mean value of the solution $S$ as a function of time, but we will use the fine grid and local problems
to construct the parameters in the proposed upscaling model. 

\subsection{Oversampling domains}
The use of oversampling domains plays a crucial role in our nonlinear upscaling technique. Two types of oversampling domains will be used.
The first type of oversampling domains, simply called oversampling domain, will be used for local problems. These local problems will be solved on oversampling regions,
and the parameters in the nonlinear upscaled model will be determined based on the these local solutions. For the transport type problems considered in this paper,
the size of the oversampling domains is determined by a suitable coarse time scale and the speed of propagation.
The second type of oversampling domains, called double oversampling domain, will be used for the purpose of imposing boundary conditions for local problems. 
This kind of domains allows an artificial layer so that artificial boundary condition imposed on them will not affect the solution in the interior. Next,
we give the definitions. 

\begin{figure}[htbp!]
  \centering
\includegraphics[width=0.6\linewidth]{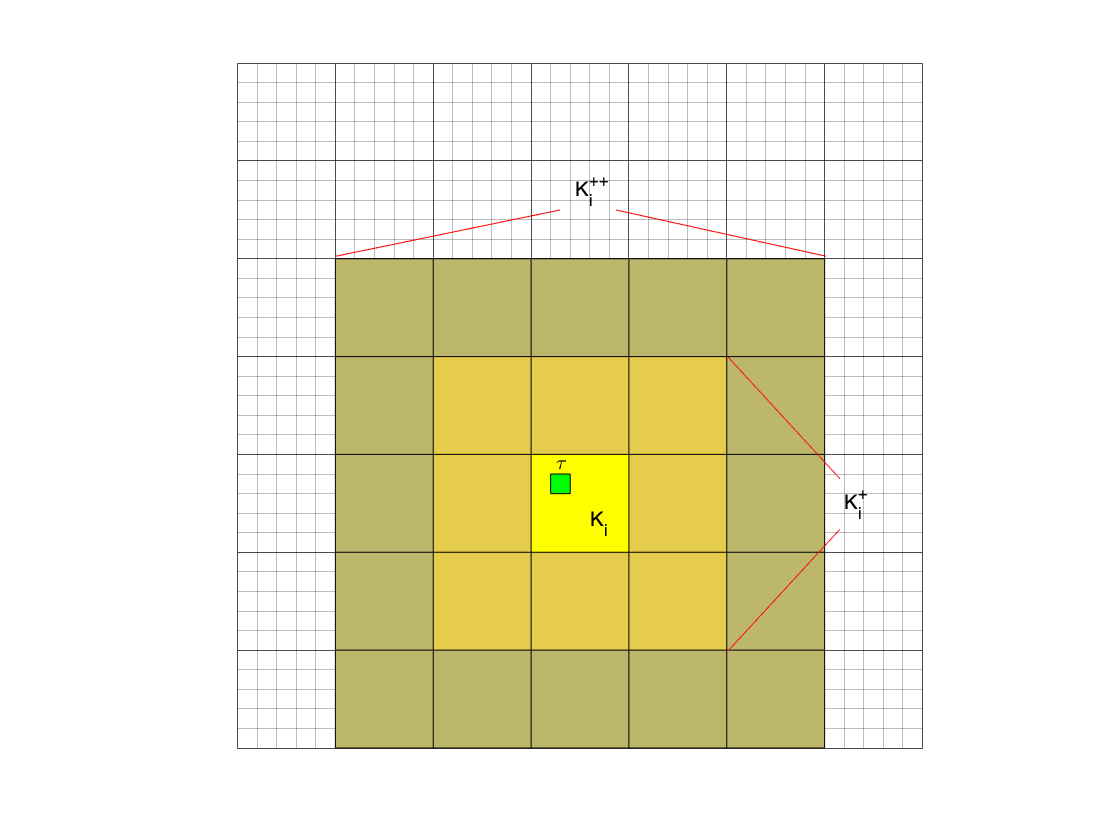}
  \caption{An illustration of the coarse grid, fine grid and oversampling domains.}
  \label{fig:testfig}
\end{figure}

 Let $K$ be a coarse element in $\mathcal{T}_H$. 
 We will construct an oversampling domain $K^+$ and a double oversampling domain $K^{++}$.
 The oversampling domain $K^+$ is obtained by enlarging $K$ by several coarse grid layers, while the double oversampling domain $K^{++}$ is obtained by enlarging $K^+$
 by several coarse grid layers. 
 An illustration of the fine grid, coarse grid, and oversampling domain are shown in Figure~\ref{fig:testfig}. 
 In this example, the oversampling domain $K^+$ is obtained by enlarging $K$ by one coarse grid layer, and the double oversampling domain $K^{++}$ is obtained by enlarging
 $K^+$ by one coarse grid layer.
 The oversampling domain $K^+$ contains totally nine coarse rectangles and double oversampling domain $K^{++}$ contains twenty-five coarse rectangles if the coarse rectangle is not near the boundary. These oversampling domains are used for the computation of local problem. 

\section{Method description}\label{sec:method}
\subsection{Algorithm}
In this section, we will present the construction of our proposed method. The main idea is based on the nonlinear NLMC framework proposed in \cite{chung2018NLNLMC}. We will apply the general concept in \cite{chung2018NLNLMC} but with some new concepts specific for the transport problem considered in this paper.  We consider the following nonlinear transport problem, in which we assume the sources $f=0$ and $g=0$ to simplify the notations. The proposed method can be easily extended to the case of nonzero sources.  We consider 
\begin{align*}
\frac{\partial S}{\partial t}+\text{div}\left(v  \lambda(S)\right)&=0, \qquad\quad~ \text{ in } (0,T)\times \Omega, \\
S&=0, \qquad \quad~\text{ on } (0,T) \times \Gamma, \\
S(0,x)&=S_0(x), \qquad \text{ in } \Omega.
\end{align*}
Our goal is to design a nonlinear upscaling method that can give the cell averages of the solution $S(t,x)$ on the coarse grid $\mathcal{T}_H$. 

For the time discretization, we will apply the backward Euler scheme. Let $\Delta t >0$ be the time step size and let $t_n = n\Delta t$. The backward Euler scheme reads
\begin{equation}
\label{eq:pde}
\gamma S^{n+1} + \text{div}\left(v  \lambda(S^{n+1})\right) = \gamma S^n
\end{equation}
where $\gamma = \frac{1}{\Delta t}$ and $S^n = S(t_n,x)$. Our proposed nonlinear upscaling method is designed to solve this equation (\ref{eq:pde}).
That is, given the mean values of the solution at time $t_n$, we will solve the mean values of the solution at $t_{n+1}$.
We use the notation $\bar{S}_{\alpha}^n$ to denote the average of $S^n$ on the coarse cell $K_{\alpha} \in \mathcal{T}_H$.
Our proposed scheme will compute $\{ \bar{S}_{\alpha}^{n+1}\}$ by using $\{ \bar{S}_{\alpha}^n\}$.

Assume that $\{ \bar{S}_{\alpha}^n\}$ are known. 
Following the general framework proposed in \cite{chung2018nlmc,zhao2020analysis,chung2018NLNLMC}, we will solve a local problem to find a downscaling function. 
To construct this downscaling function, we consider a set of values $\{ \bar{S}_{\beta}\}$, where each $\bar{S}_{\beta}$ represents the average value of a certain function on 
the coarse cell $K_{\beta}$. We remark that these values $\bar{S}_{\beta}$ will be the required values $\{ \bar{S}^{n+1}_{\beta}\}$ after solving our upscaled model. 
To get the required downscaling function, we will solve a local problem to find a local downscaling function $\tilde{\psi}$ 
such that the mean values of $\tilde{\psi}$ on the coarse cells are given by the values $\{ \bar{S}_{\beta}\}$. 
Specifically, 
let $K^{++}_\alpha$ be a double oversampling domain of the coarse cell $K_\alpha\in \mathcal{T}_H$. 
We consider the local problem:
\begin{align}
\gamma \tilde{\psi}+\text{div}\left(v \lambda(\tilde{\psi})\right)&=\gamma \bar{S}^n, \qquad \text{in $K^{++}_\alpha$}, \label{eq:loc1} \\
 \tilde{\psi}&=0, \qquad \text{on $\partial^- K^{++}_\alpha$}, \label{eq:loc2} \\
 \frac{1}{|K_\beta|}\int_{K_\beta}  \tilde{\psi} &= \bar{S}_{\beta}, \qquad \text{$K_\beta\subset K^{++}_\alpha$ and $K_\beta\in \mathcal{T}_H$} \label{eq:loc3}
\end{align}
where $\partial^- K^{++}_\alpha$ is the inflow boundary of $K^{++}_{\alpha}$. 
The above problem is solved numerically on the fine grid $\mathcal{T}_h$. The solution $\tilde{\psi}$
will give the required local downscaling function. Next, 
we restrict the solution $\tilde{\psi}$ on $K^+_{\alpha}$ and 
 denote it as $\psi_\alpha$, which is the required local downscaling function whose support is $K^+_{\alpha}$.

Using these local solutions $\psi_\alpha$, we define a global downscaling function $\psi = \sum \chi_\alpha \psi_\alpha$,
where $\{ \chi_{\alpha}\}$ is a set pf partition of unity functions corresponding to the overlapping partition $\{ K_{\alpha}^+\}$ of the domain. 
We remark that this global function $\psi$ depends on the mean values $\{ \bar{S}_\beta \}$.
Now we solve (\ref{eq:pde}) using this global downscaling function. We will integrate the equation (\ref{eq:pde}) on each coarse element $K_i$
and replace the true solution $S^{n+1}$ by $\psi$. In particular, we have
\begin{equation}
\label{eq:scheme}
\gamma \int_{K_i} \psi + \int_{\partial K_i} \lambda(\psi) v\cdot n = \gamma |K_i| \bar{S}_i^n.
\end{equation}
Here, we use the notation $\bar{S}^n$ to denote the values $(\bar{S}^n_1, \cdots, \bar{S}^n_{N_K})$.
We will solve this nonlinear equation and obtain the unknowns $\{ \bar{S}_{\beta}\}$. These values will be the required solution at the time step $t_{n+1}$.
That is, we set $\bar{S}^{n+1}_{\beta}$ to be the solution obtained by solving (\ref{eq:scheme}).

To numerically solve (\ref{eq:scheme}), we apply a fixed point type iteration. 
Let $\bar{S}^{(m)}=(\bar{S}^{(m)}_1,\cdots,\bar{S}^{(m)}_{N_K})$ be the $m$-th iterate. 
We can then construct the global downscaling function $\psi(\bar{S}^{(m)}_{1},\cdots,\bar{S}^{(m)}_{N_K})$  using the values $\bar{S}^{(m)}_{i}$.
Finally we perform the following updating procedure
to find the mean values $\bar{S}^{(m+1)}_{i}$:
\begin{align}
\begin{cases}
F^{(m)}_{i}(\overline{S}^{(m)})=\gamma \int_{K_i} \psi(\overline{S}^{(m)})- \gamma |K_i| \bar{S}^n_i+\int_{\partial K_i}  \lambda(\psi(\overline{S}^{(m)}) v\cdot n\\
\bar{S}^{(m+1)}=\bar{S}^{(m)}-F^{(m)}
\end{cases}
\label{eq:fixedpt}
\end{align}
with a given
initial guess $\bar{S}^{(0)}$. Assume that the above iteration converges at the $M$-th iteration. 
We will set $\bar{S}^{n+1} = \bar{S}^{(M)}$.
We remark that each iteration of the above updates requires the construction of the global downscaling function $\psi(\bar{S}^{(m)}_{1},\cdots,\bar{S}^{(m)}_{N_K})$,
which requires the solution of many local nonlinear problems (\ref{eq:loc1})-(\ref{eq:loc3}). 
With the use of the proposal deep learning approach, which will be presented in Section~\ref{sec:deep}, this step becomes much more efficient. 

\subsection{Computation of local problems}

We will use the Newton's method to solve the local problem (\ref{eq:loc1})-(\ref{eq:loc3}).
We recall that the value $\bar{S}^n$ is fixed. For a given $\{ \bar{S}_{\beta}\}$, the problem (\ref{eq:loc1})-(\ref{eq:loc3}) gives
a function $\tilde{\psi}$ defined on $K_{\alpha}^{++}$. 
We note that the value $\{ \bar{S}_{\beta}\}$ can be considered as a piecewise constant function, denoted as $\bar{S}$, defined on the coarse grid supported in $K_{\alpha}^{++}$.
Thus, the system (\ref{eq:loc1})-(\ref{eq:loc3})
defines a map which takes a piecewise constant function $\{ \bar{S}_{\beta}\}$ as input
and returns a fine scale function $\tilde{\psi}$. In our numerical simulations, we take $\tilde{\psi}$ to be piecewise constant function on the fine grid. 
In order to ensure the constraint (\ref{eq:loc3}) is satisfied, we will introduce an additional variable $\mu$.
This variable $\mu$ is a piecewise constant function on the coarse grid whose support is $K_{\alpha}^{++}$.
In addition, the system (\ref{eq:loc1})-(\ref{eq:loc3}) will be solved on the fine grid by
using the upwind finite volume scheme. In the following, we will present the mathematical details. 

The equation (\ref{eq:loc1}) will be discretized by the upwind finite volume scheme. 
That is, for each fine grid element $\tau_j \in\mathcal{T}_h$, we will apply the upwind finite volume scheme to (\ref{eq:loc1}). 
This motivates us to define the following two operators:
\begin{align*}
\begin{cases}
F_{\text{eqn},j}(\tilde{\psi},\mu)=\gamma |\tau_j| \tilde{\psi}_j+\int_{\partial \tau_j} \lambda(\tilde{\psi}_{\text{up}}) v\cdot n -\int_{\tau_j} \bar{S}^n -\int_{\tau_j} \mu,
                                                                                                                                                                                                                                                                                                                                                                                                                                                                                                                                                                                                                                                                                                                                                                                                                                                                                                                                                                                                                                                                                                                                                                                                                                                                                                                                                                                                                                                                                                                                                                                                       &~ \forall \tau_j\subset K^{++}_\alpha, \\
F_{\text{mean},\beta}(\tilde{\psi},\mu)=- |\tau_j|\sum_{\tau_j\subset K_\beta}\tilde{\psi}_j+ |K_{\beta}|\bar{S}^{(m)}_\beta, &~  \forall K_\beta\subset K^+_\alpha,
\end{cases}
\end{align*}
where $\tilde{\psi}_{\text{up}}$ denotes the upwind flux. We use the notation $F_{\text{eqn}}$ to denote the vector whose $j$-th component is $F_{\text{eqn},j}$, 
and use the notation $F_{\text{mean}}$ to denote the vector whose $\beta$-th component is $F_{\text{mean},\beta}$.
We remark that $F_{\text{eqn},j}$ gives the discretization of (\ref{eq:loc1}) on the fine grid cell $\tau_j \in\mathcal{T}_h$,
and $F_{\text{mean},\beta}$ corresponds to the equation (\ref{eq:loc3}) on the coarse element $K_{\beta}$.
Next, we define
 \begin{align*}
F(\tilde{\psi},\mu)=\begin{cases}
F_\text{eqn}(\tilde{\psi},\mu)\\
F_\text{mean}(\tilde{\psi},\mu).
\end{cases}
\end{align*}
The goal is to apply the Newton's method to solve $F(\tilde{\psi},\mu)=0$, which gives (\ref{eq:loc1}) and (\ref{eq:loc3}).
The boundary condition (\ref{eq:loc2}) is imposed by the numerical flux $\tilde{\psi}_{\text{up}}$.

In order to apply the Newton's method, we need
the Jacobian matrix $J$ of $F$, which is given by 
 \begin{align*}
 \displaystyle
J(\tilde{\psi},\mu)&=\begin{bmatrix} \frac{\partial F_\text{eqn}}{\partial \tilde{\psi}} &\frac{\partial F_\text{eqn}}{\partial \mu}\\
 \frac{\partial F_\text{mean}}{\partial \tilde{\psi}} &\frac{\partial F_\text{mean}}{\partial \mu}
\end{bmatrix}\\
&=\begin{bmatrix} \frac{\partial F_\text{eqn}}{\partial \tilde{\psi}} & B\\
B^T&0
\end{bmatrix}
\end{align*}
where
 \begin{align*}
 \frac{\partial F_{\text{eqn},j}}{\partial \tilde{\psi}_i} =\gamma |\tau_j| \delta_{i,j}+\int_{\partial \tau_j} \lambda'(\tilde{\psi}_{\text{up}}) \, v\cdot n \, \delta_{i,j} \quad \text{and} \quad  B_{\beta,j}=-|\tau_j|\delta^K_{\beta,j}
 \end{align*}
where $\delta_{i,j}$ is the delta function and  
 \begin{align*}
 \delta^K_{\beta,j}=
 \begin{cases}
 1 \quad \text{if } \tau_j\subset K_\beta\\
  0 \quad \text{otherwise}
 \end{cases}
 \end{align*}
The
Newton's method is then given by 
 \begin{align}
 \begin{pmatrix}  \tilde{\psi}^{(n+1)}\\ \mu^{(n+1)} \end{pmatrix}&=\begin{pmatrix}  \tilde{\psi}^{(n)}\\ \mu^{(n)}\end{pmatrix}-J(\tilde{\psi}^{(n)},\mu^{(n)})^{-1} F(\tilde{\psi}^{(n)},\mu^{(n)})
 \label{eq:newtonF}
 \end{align}
 with a suitable initial guess $(\tilde{\psi}^{(0)},\mu^{(0)})$. 
We will stop the iteration if the value $F(\tilde{\psi}^{(n)},\mu^{(n)})$ is sufficiently small. 
This procedure gives the solution of (\ref{eq:loc1})-(\ref{eq:loc3}).

 \subsection{Deep neural network model for local downscaling functions}\label{sec:deep}
We let the function $\mathcal{N}$ be a network of $L$ layers, $x$ be the input and $y$ be the
corresponding output. We write
\begin{align*}
\mathcal{N}(x ; \theta)=\sigma\left(W_{L} \sigma\left(\cdots \sigma\left(W_{2} \sigma\left(W_{1} x+b_{1}\right)+b_{2}\right) \cdots\right)+b_{L}\right)
\end{align*}
where $\theta :=\left(W_{1}, \cdots, W_{L}, b_{1}, \cdots, b_{L}\right)$, $W_i$ are the weight matrices and $b_i$ are the bias vectors, and $\sigma$ is the activation function. A neural network describes the
connection of a collection of nodes (neurons) sit in successive layers. The output neurons in each layer are
simultaneously the input neurons in the next layer. The data propagate from the input layer to the output
layer through hidden layers. The neurons can be switched on or off as the input is propagated forward
through the network.

Suppose we are given a collection of sample pairs $\{ \left(x_{j}, y_{j}\right) \}$. The goal is then to find $\theta^{*}$ by solving an optimization problem 
\begin{align*}
\theta^{*}=\underset{\theta}{\operatorname{argmin}} \frac{1}{N_s} \sum_{j=1}^{N_s}\left\|y_{j}-\mathcal{N}\left(x_{j} ; \theta\right)\right\|_{2}^{2}
\end{align*}
where $N_s$ is the number of the samples. Here, the function $\frac{1}{N_s} \sum_{j=1}^{N_s}\left\|\mathbf{y}_{j}-\mathcal{N}\left(\mathbf{x}_{j} ; \theta\right)\right\|_{2}^{2}$ is known as the loss function. One needs to select suitable number of layers, number of neurons in each layer, the activation function, the loss function and the optimizers for the network. 

We will use a deep neural network $\mathcal{N}$ to model the process of constructing downscaling functions. 
Recall that the local downscaling function is defined by (\ref{eq:loc1})-(\ref{eq:loc3}).
We note that there are two sources of inputs. The first input contains the values $\{ \bar{S}^n_i\}$ from the time step $n$. 
The second input contains the values $\{ \bar{S}_{\beta}^{(m)}\}$ in the fixed point iteration (\ref{eq:fixedpt}).
Moreover, we observe that the outputs that we need are the values of the global downscaling function $\psi$ 
restricted to the coarse grid edges, since only these values are used in the numerical scheme (\ref{eq:scheme}).

Therefore, we use
the following choices in our deep neural network:
\begin{itemize}
\item Input: $x=\{\bar{S}^{(m)},\bar{S}^n \}$. 
\item Output: $y=\psi|_{\cup_j E_j}$, where $E_j$ is the $j$-th coarse edge. This is the restriction of $\psi$ on all coarse edges. 
\item Sample pairs: $N_s=5000$ sample pairs of $\left(x_{j}, y_{j}\right)$ are collected.
\item Standard loss function: $\frac{1}{N_s} \sum_{j=1}^{N_s}\left\|y_{j}-\mathcal{N}\left(x_{j} ; \theta\right)\right\|_{2}^{2}$.
\item Activation function: The popular ReLU function (the rectified linear unit activation function) is a
common choice for activation function in training deep neural network architectures.
\end{itemize}

As for the input of the network, we use $x=\{\bar{S}^{(m)},\bar{S}^n \},$ which are the vectors containing the approximate mean value $\bar{S}^{(m)}$ from the current iteration in (\ref{eq:fixedpt})
and the mean value $\bar{S}^n$ at the $n$-th time step of the scheme (\ref{eq:scheme}). The input $x$ is a random vector such that each entry is
ranged from 0 to 1. The range is based on the range of the initial condition. Since the range of the initial condition will affect the range of the mean value of the next time step, we choose the range of random input vector $x$ to be $[0,1]$. 

The corresponding output data are $y=\psi|_{\cup_j E_j},$ which contains the values of the global downscaling function on each coarse edge. 
In order to obtain the data $y$, 
we use randomly generated input $x$ and solve the local system (\ref{eq:loc1})-(\ref{eq:loc3}) together with $\psi = \sum \chi_{\alpha} \psi_{\alpha}$
to obtain the global downscaling function $\psi$, which will give the output data $y$.

In between the input and output layer, we test on $2$ hidden layers with the number of neurons ranged between $4N_K$ and $8N_K$
in each hidden layer,  where $N_K$ is the number of the coarse elements. We note that the size of input vector $x$ is $2N_K$ and the size of output vector $y$ is $n_E N_E$, where $N_E$ is the number of coarse grid edges and $n_E$ is the number of fine grid cells on the coarse edge $E$. In the training, there are $N_s=5000$ sample pairs of $\left(\mathbf{x}_{j}, \mathbf{y}_{j}\right)$ collected.

In between layers, we need the activation function. The ReLU function (rectified linear unit activation function) is a popular choice for activation function in training deep neural network architectures. And we will use the standard loss function: $\frac{1}{N_s} \sum_{j=1}^{N_s}\left\|y_{j}-\mathcal{N}\left(x_{j} ; \theta\right)\right\|_{2}^{2}$. As for the training optimizer, we use AdaMax, which is a stochastic gradient descent (SGD) type algorithm well-suited for high-dimensional parameter space, in minimizing the loss function.


\section{Numerical examples}\label{sec:num}

In this section, we will present some numerical examples to show the performance of our proposed deep learning based nonlinear upscaling method.
In our simulations, we will take $\lambda(S) = S^2$. 
Moreover, the oversampling region $K^+$ is chosen by enlarging the coarse cell $K$ by one coarse grid layer,
and the double oversampling region $K^{++}$ is chosen by enlarging $K^+$ by one coarse grid layer.
This choice of oversampling regions is motivated by the coarse time step size and the velocity of propagation. 
We will show the performance by using our upscaling method without using deep learning in the construction of the local downscaling functions.
In this case, the error will only come from the approximation of the upscaling method. In addition, we will present the performance of our method with the use of deep neural network to approximate
the local downscaling functions. In this case, the error in this deep neural network approximation will be added to the solutions. However, the computational efficiency improves significantly. 
In both cases, we observe that our proposed upscaling method is able to give accurate numerical approximations. 

\subsection{Example 1}
\begin{figure}
\label{fig:f1}
\centering
  \includegraphics[scale=0.3]{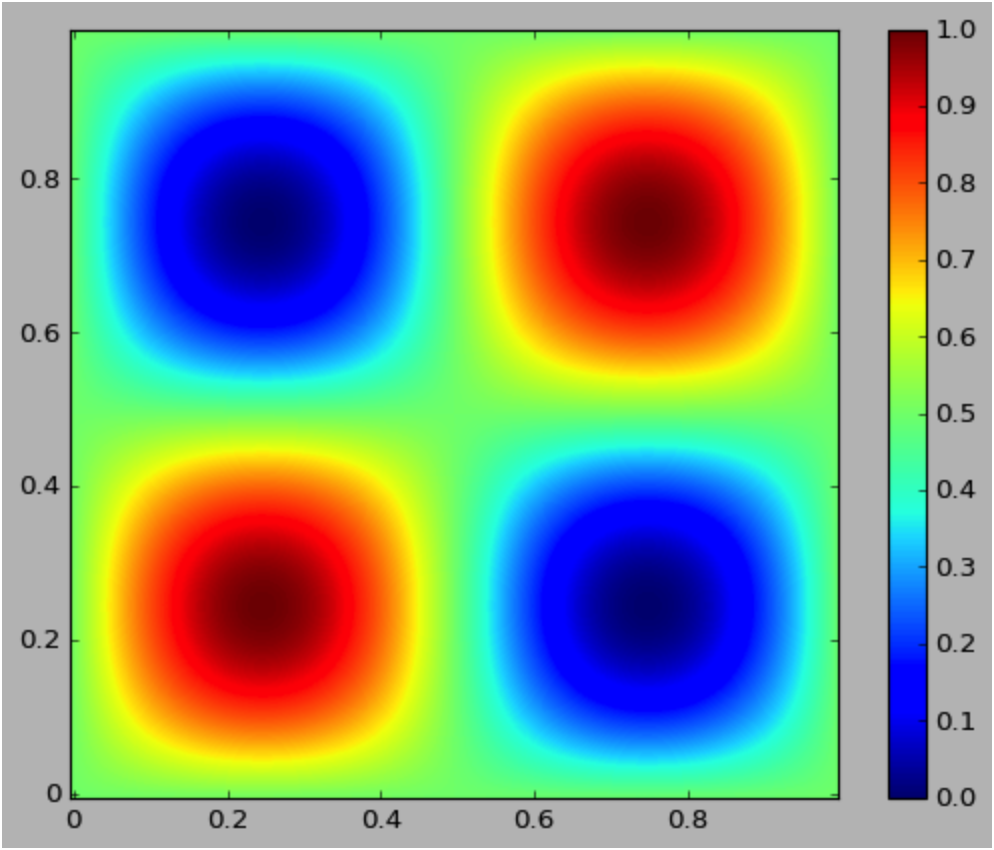}
\caption{Initial condition for Example 1.}
\end{figure}

For the first numerical example, we will take the time step as $\Delta t=0.1$ and use the constant velocity  $v=(1,1)$ to test our algorithm. As a result, the constant $\gamma=\frac{1}{\Delta t}=10$. The coarse mesh is $20\times 20$ and the fine mesh is $100\times 100$.  
In this case, due to the finite speed of propagation, the size of the oversampling region $K^+$ is enough to capture the solution at the next time step
 from the solution at the previous time step originated in $K$.  
Also, we choose the initial condition as $S_0(x,y)=(1+\sin(2\pi x)\sin(2\pi y))/2$. Figure~\ref{fig:f1} shows a plot of  $S_0$. The initial condition is bounded between $0$ and $1$. For the numerical computation of equation (\ref{eq:scheme}) to find the mean value $\bar{S}^{n+1}$ for the next time step, we use the fixed point type iteration (\ref{eq:fixedpt}) so we need a stopping criterion to stop the iteration. In our simulations, we choose the condition $\|\bar{S}^{(m+1)}-\bar{S}^{(m)}\|\leq 0.0001$ as the stopping criterion.

\begin{figure}
\centering
  \includegraphics[scale=0.4]{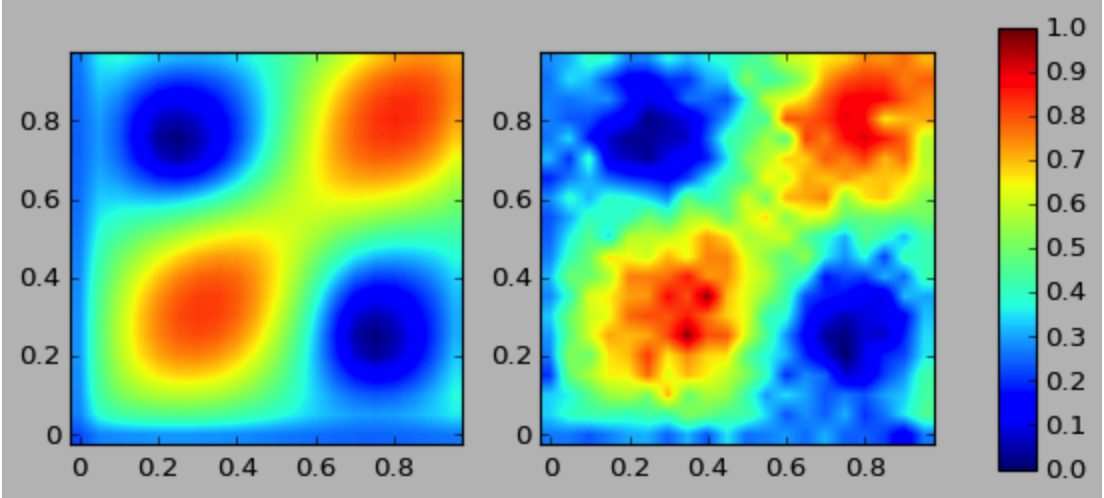}
    \includegraphics[scale=0.4]{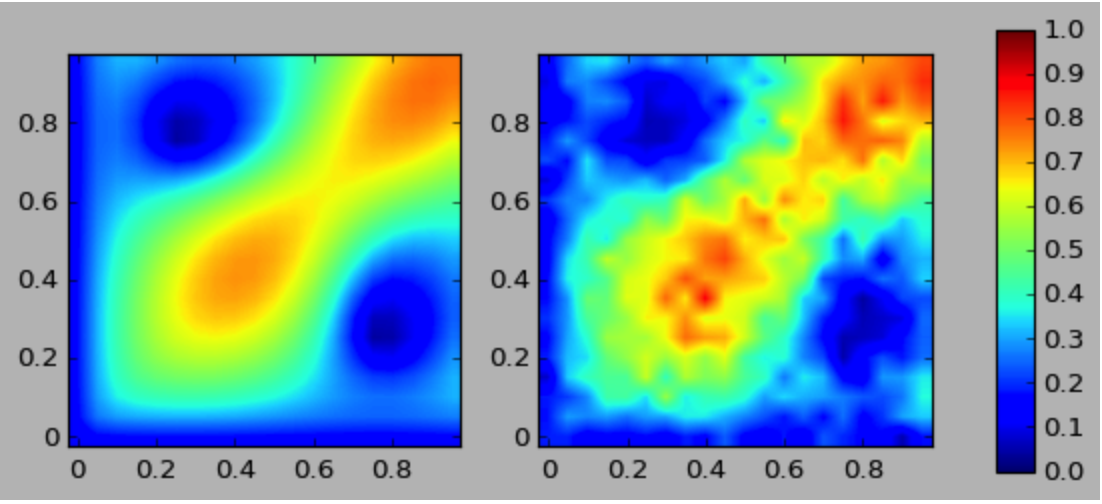}
      \includegraphics[scale=0.4]{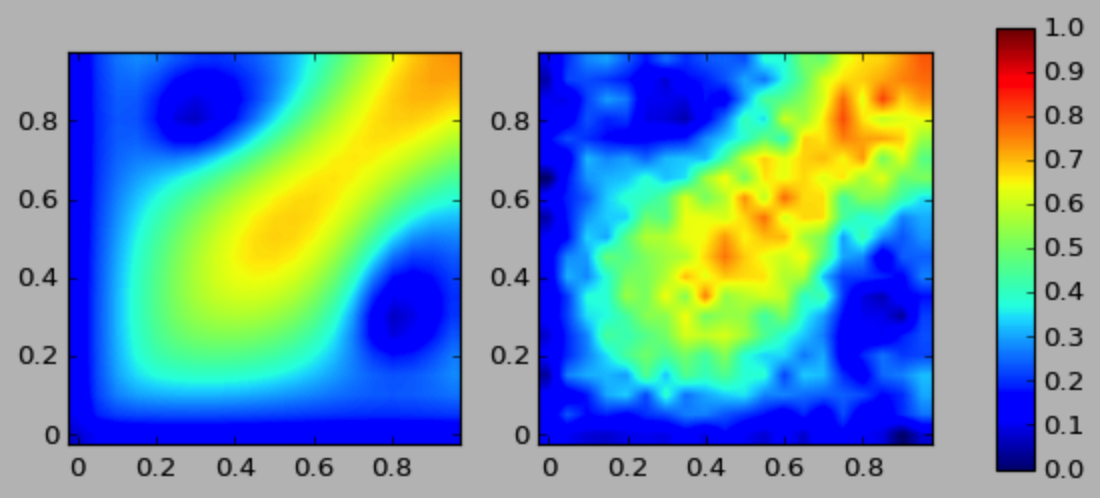}
        \includegraphics[scale=0.4]{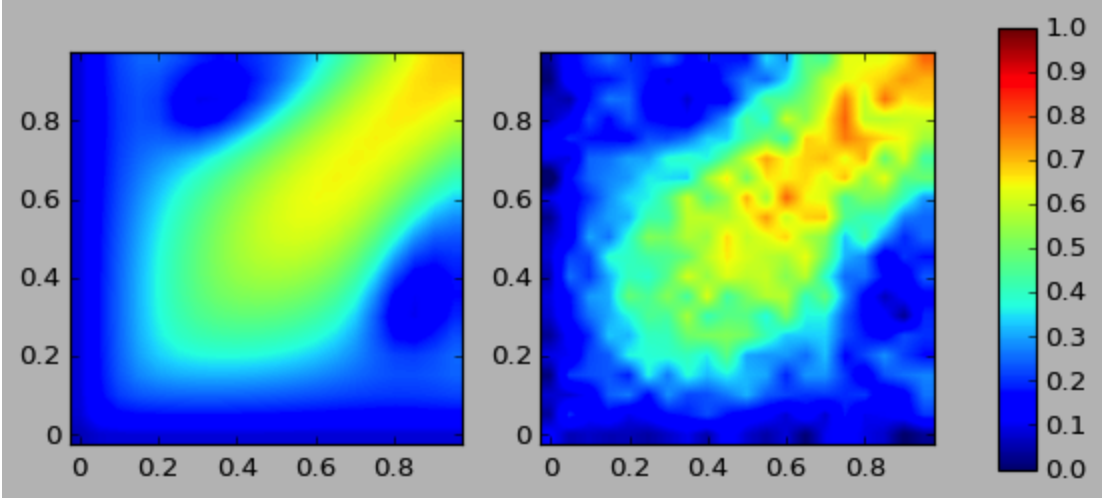}
          \includegraphics[scale=0.4]{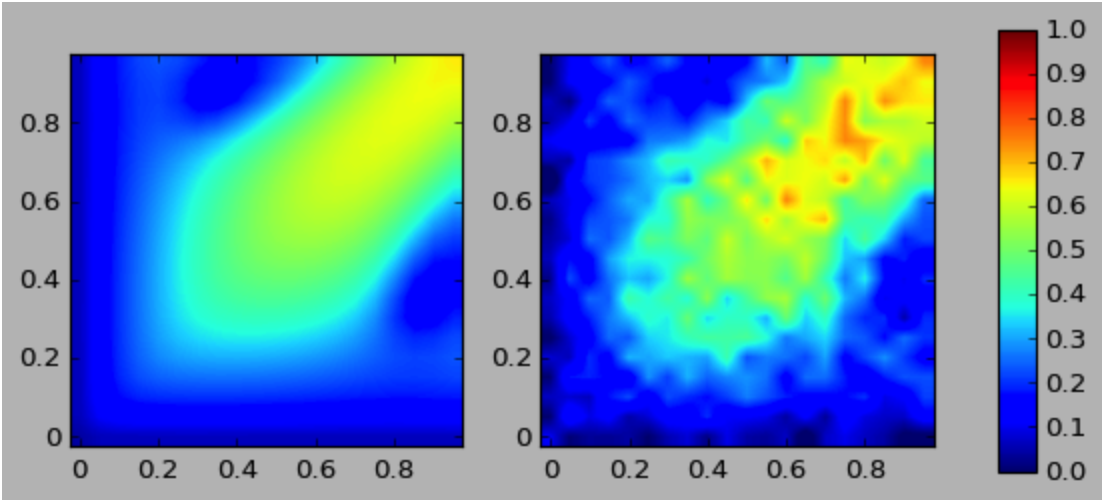}
\caption{Solutions at the first 5 time steps $t=0.1,0.2, \cdots, 0.5$ for Example 1. Left: upscaled reference solution. Right: numerical solution using deep learning based upscaling.}
  \label{fig:solutionpic1}
\end{figure}

In this numerical example, we test the algorithm by running 5 time steps.  Figure~\ref{fig:solutionpic1} shows the solutions for the first 5 time steps of our proposed scheme using neural network model. 
We also show the average value of the reference solution on the coarse grid in this figure. From these results, 
we observe a very good agreement between these two solutions. 
We have also observed that the numerical solution is bounded between $0$ and $1$, which holds for the reference solution.
We remark that we did not impose any bound preserving properties in the training of our deep neural network.
An improved network with this desirable property will be developed in a future work. 

%

\begin {table}[H]
 \caption{Root mean square error for Example 1 using neural network model}
  \label{fig:table1.1}
\begin{center}
 \begin{tabular}{|c|c|} 
 \hline
time step  & relative error of the mean value \\ [0.5ex] 
 \hline\hline
$t=0.1$&0.04379  \\ 
 \hline
$t=0.2$&0.04728  \\ 
 \hline
$t=0.3$& 0.05264  \\ 
 \hline
$t=0.4$&0.05837  \\ 
 \hline
 $t=0.5$&0.06569  \\ 
 \hline
\end{tabular}
\end{center}
\end {table}

For the proposed neural network model, the relative error of the mean value at the first time step is 0.04379. The relative error of the mean value of at the fifth time step is 0.06569. Table~\ref{fig:table1.1} shows the relative errors for the mean value of the solution using our neural network model. From the results in this table, we observe that
our scheme gives a very good performance. 
We also observe that the relative error of the mean value rises slowly when the time increases. 
In terms of the computational efficiency, 
the time used is about 45s for computation of one time step for our neural network model. 
As a comparison, if we compute the downscaling functions without the use of neural networks, we need about 4500s for computation of one time step. Most of the computational time 
is consumed on solving the local downscaling functions. So, we see that the use of deep neural network model can reduce the computational time without sacrificing the accuracy.

\begin {table}[H]
 \caption{Root mean square error for Example 1 using the upscaled model without neural network}
  \label{fig:table1.2}
\begin{center}
 \begin{tabular}{|c|c| c|} 
 \hline
time step  & relative error of the mean value \\ [0.5ex] 
 \hline\hline
$t=0.1$ & 0.03210  \\ 
 \hline
$t=0.2$ & 0.04070  \\ 
 \hline
$t=0.3$ & 0.04574  \\ 
 \hline
$t=0.4$ & 0.05227  \\ 
 \hline
 $t=0.5$ & 0.05887  \\ 
 \hline
\end{tabular}
\end{center}
\end {table}

We also report the results using the downscaling functions without deep neural network in the construction of the proposed upscaling model. The relative error of  the mean value at the 1st time step is 0.03210. The relative error of the mean value at the 5th time step is 0.05887. Table~\ref{fig:table1.2} shows the relative errors of the mean value when our method is applied without the use of deep neural network. The performance of solving the local downscaling functions numerically without the use of deep neural network and solving the local downscaling functions using neural network model are similar
as we can see from Table~\ref{fig:table1.1} and Table~\ref{fig:table1.2}. 

\begin {table}[H]
 \caption{Root mean square error for Example 1 using the finite volume method}
  \label{fig:table1.3}
\begin{center}
 \begin{tabular}{|c|c|} 
 \hline
time step  & relative error of the mean value \\ [0.5ex] 
 \hline\hline
$t=0.1$ & 0.06180  \\ 
 \hline
$t=0.2$ & 0.08007  \\ 
 \hline
$t=0.3$ & 0.09705  \\ 
 \hline
$t=0.4$ & 0.10891 \\ 
 \hline
 $t=0.5$ & 0.11433  \\ 
 \hline
\end{tabular}
\end{center}
\end {table}

We next compare our scheme with
a standard coarse-grid finite volume method with upwind flux. The relative error of the mean value at the 1st time step is 0.06180. The relative error of the mean value at the 5th time step is 0.11433. Table~\ref{fig:table1.3} shows the  relative errors of the mean value when a standard upwind finite volume method is applied on the coarse grid. We notice that the performance of our method with the use of deep neural network model is better than a simple application of a coarse-grid finite volume method.

\subsection{Example 2}
\begin{figure}[ht]
   \centering
  \includegraphics[width= 0.7\linewidth]{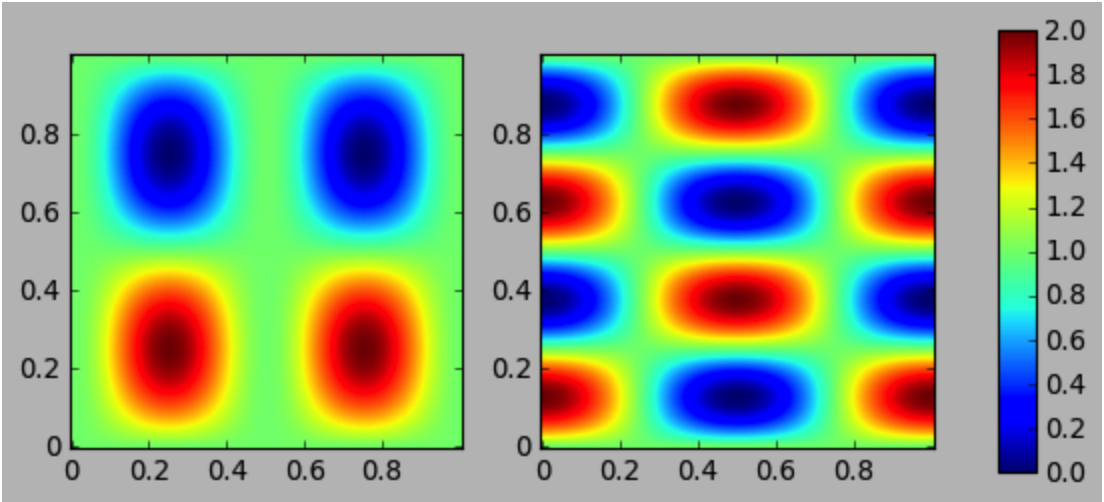}
  \caption{The velocity $v(x,y) = (v_1,v_2)$ for Example 2. The left figure is $v_1$ and the right figure is $v_2$.}
    \label{fig:v2}
\end{figure}

For the second numerical example, we will also take the time step as $\Delta t=0.1$ and use a velocity given by
$$v(x,y)=( 1+\sin^2(2\pi x)\sin(2\pi y),1+\sin(4\pi x)\cos(2\pi y))$$
 to test our algorithm. 
Figure~\ref{fig:v2} shows the plot of  $v(x,y)$.
 The coarse mesh is also chosen as $20\times 20$ and the fine mesh is $100\times 100$.  Since the velocity is bounded between 0 and 2, the oversampling region $K^+$ is able to keep all information originated from the coarse element $K$ within one coarse time step. 
 Also, we choose the same initial condition as $S_0(x,y)=(1+\sin(2\pi x)\sin(2\pi y))/2$. 
 Same as Example 1, we choose the condition $\|\bar{S}^{(m+1)}-\bar{S}^{(m)}\|\leq 0.0001$ as the stopping criterion.

\begin{figure}
\centering
  \includegraphics[scale=0.4]{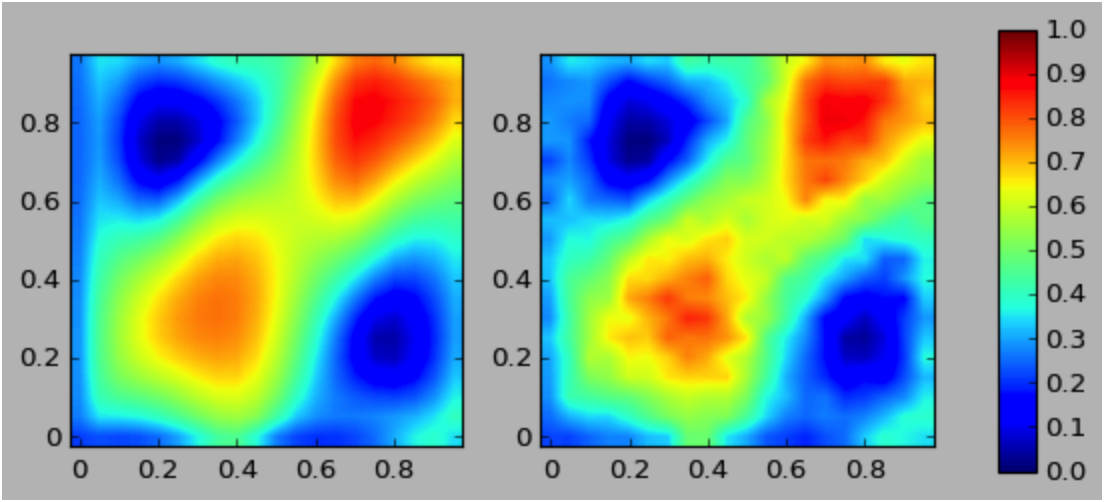}
    \includegraphics[scale=0.4]{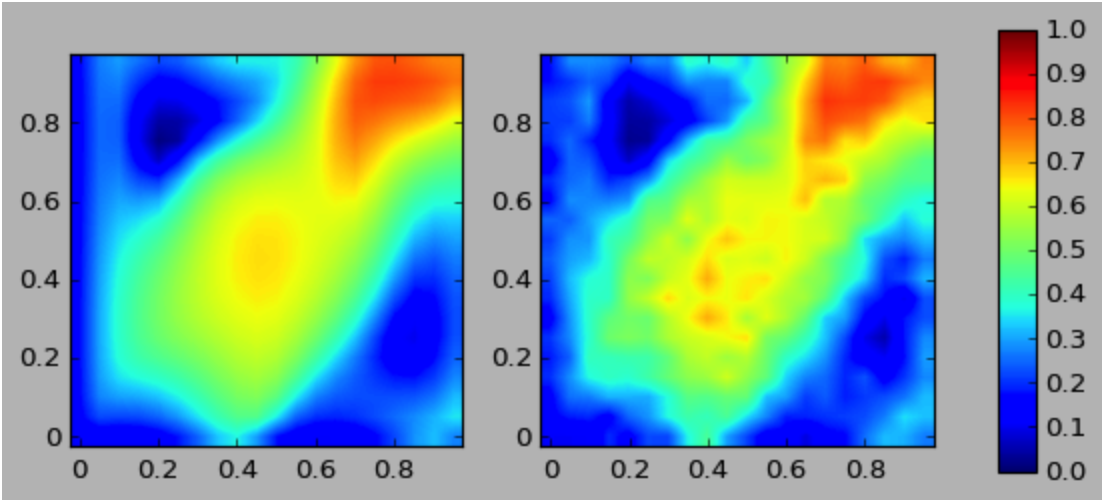}
      \includegraphics[scale=0.4]{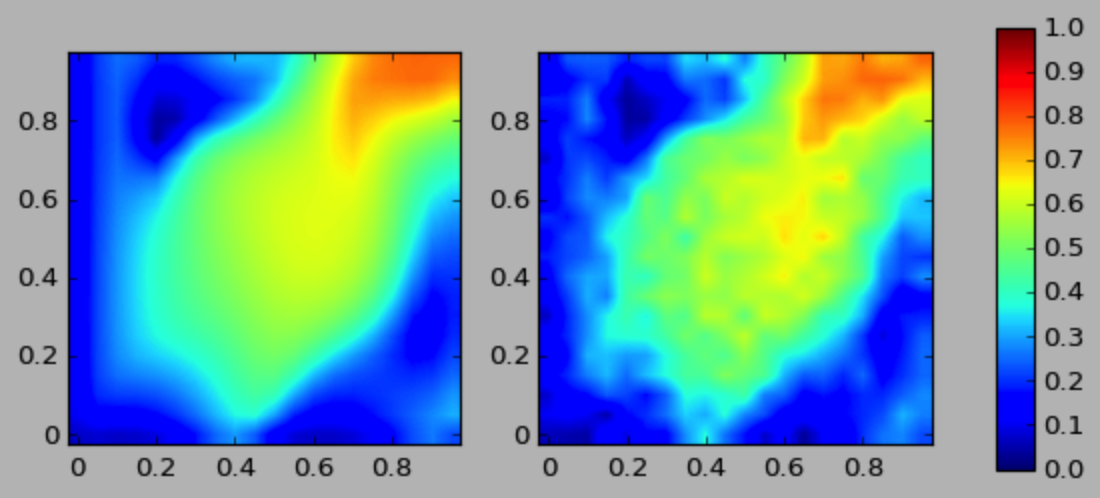}
        \includegraphics[scale=0.4]{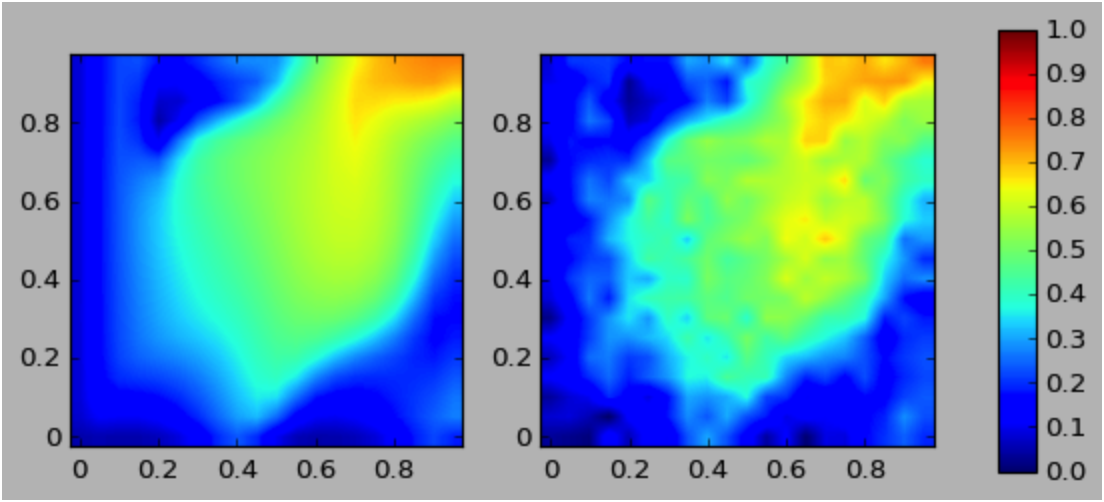}
          \includegraphics[scale=0.4]{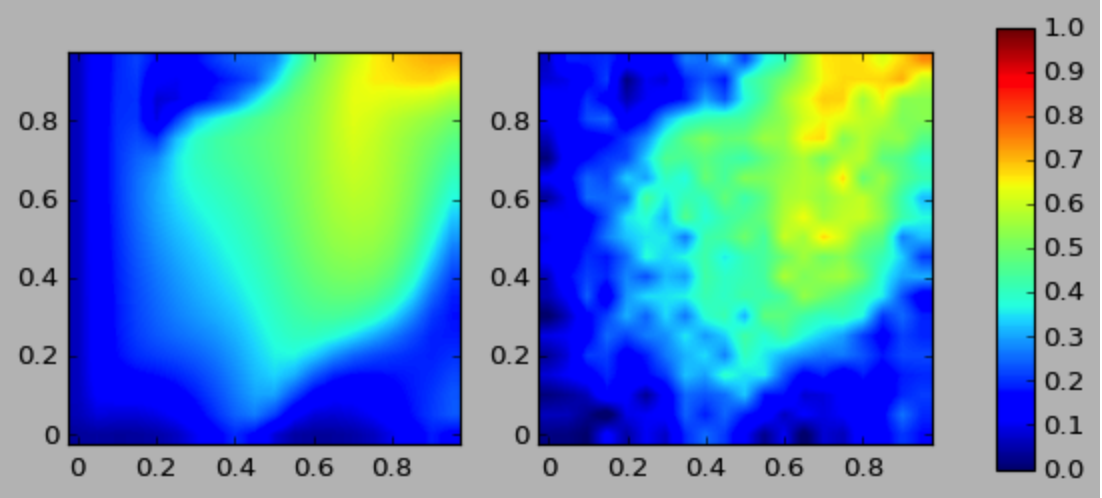}
\caption{Solutions at the first 5 time steps $t=0.1,0.2, \cdots, 0.5$ for Example 2. Left: upscaled reference solution. Right: numerical solution using deep learning based upscaling.}
\label{fig:solutionpic2}
\end{figure}

In this numerical example, we test the algorithm with 5 time steps.  Figure~\ref{fig:solutionpic2} shows the numerical solutions for
the first 5 time steps computed using our proposed scheme with deep neural network model. 
In the same figure, we also show the corresponding average values of the reference solution on the coarse grid. 
From these figures, we observe very good agreement between the numerical solution and the average reference solution. 

%

\begin {table}[H]
 \caption{Root mean square error for Example 2 using neural network model}
  \label{fig:table2.1}
\begin{center}
 \begin{tabular}{|c|c|} 
 \hline
time step  & relative error of the mean value \\ [0.5ex] 
 \hline\hline
$t=0.1$&0.04468  \\ 
 \hline
$t=0.2$&0.04587  \\ 
 \hline
$t=0.3$&0.04727  \\ 
 \hline
$t=0.4$&0.04795  \\ 
 \hline
 $t=0.5$&0.05655  \\ 
 \hline
\end{tabular}
\end{center}
\end {table}

For our proposed scheme with the use of deep neural network model, the relative error of the numerical solution at the first time step is 0.04468. The relative error of the numerical solution at the fifth time step is 0.05655. Table~\ref{fig:table2.1} shows the relative errors of the numerical solution computed using our proposed scheme with neural network model, and we 
observe a very good performance. 
Also the computational time is about 50s for the computation of one time step of using the proposed neural network model. As a comparison, if we compute the downscaling functions without neural networks, we need about 5500s for computation of one time step. 

\begin {table}[H]
 \caption{Root mean square error for Example 2 using the upscaled model without neural network}
  \label{fig:table2.2}
\begin{center}
 \begin{tabular}{|c|c|} 
 \hline
time step  & relative error of the mean value \\ [0.5ex] 
 \hline\hline
$t=0.1$ & 0.02904  \\ 
 \hline
$t=0.2$ & 0.03579  \\ 
 \hline
$t=0.3$& 0.04223  \\ 
 \hline
$t=0.4$ & 0.04829  \\ 
 \hline
 $t=0.5$ & 0.05499  \\ 
 \hline
\end{tabular}
\end{center}
\end {table}

As a comparison, we present the results for our scheme without neural network. The relative error of the numerical solution at the 1st time step is 0.02904. The relative error of the numerical solution at the 5th time step is 0.05499. Table~\ref{fig:table2.2} shows the relative errors of the numerical solution computed using our scheme without neural network. That is, we solve the local problem (\ref{eq:loc1})-(\ref{eq:loc3}) directly. 
We observe that the performance of using neural network and without the use of neural network are similar. 
This shows that our proposed deep neural network can increase the efficiency of the method without losing accuracy. 

\begin {table}[H]
 \caption{Root mean square error for Example 2 using the finite volume method}
  \label{fig:table2.3}
\begin{center}
 \begin{tabular}{|c|c|} 
 \hline
time step  & relative error of the mean value \\ [0.5ex] 
 \hline\hline
$t=0.1$ & 0.06116  \\ 
 \hline
$t=0.2$& 0.07144  \\ 
 \hline
$t=0.3$& 0.08546  \\ 
 \hline
$t=0.4$& 0.09330  \\ 
 \hline
 $t=0.5$ & 0.09634  \\ 
 \hline
\end{tabular}
\end{center}
\end {table}

Again, we show the performance of using a standard coarse-grid finite volume method with upwind flux. The relative error of the numerical solution at the 1st time step is 0.06116. The relative error of the numerical solution at the 5th time step is 0.09634. Table~\ref{fig:table2.3} shows the relative errors of the numerical solution computed using finite volume method. 
Comparing the results in Table~\ref{fig:table2.3} and Table~\ref{fig:table2.1},
we observe that our proposed nonlinear upscaling is able to give better numerical solutions. 

\subsection{Example 3}

\begin{figure}[H]
\centering
  \includegraphics[scale=0.3]{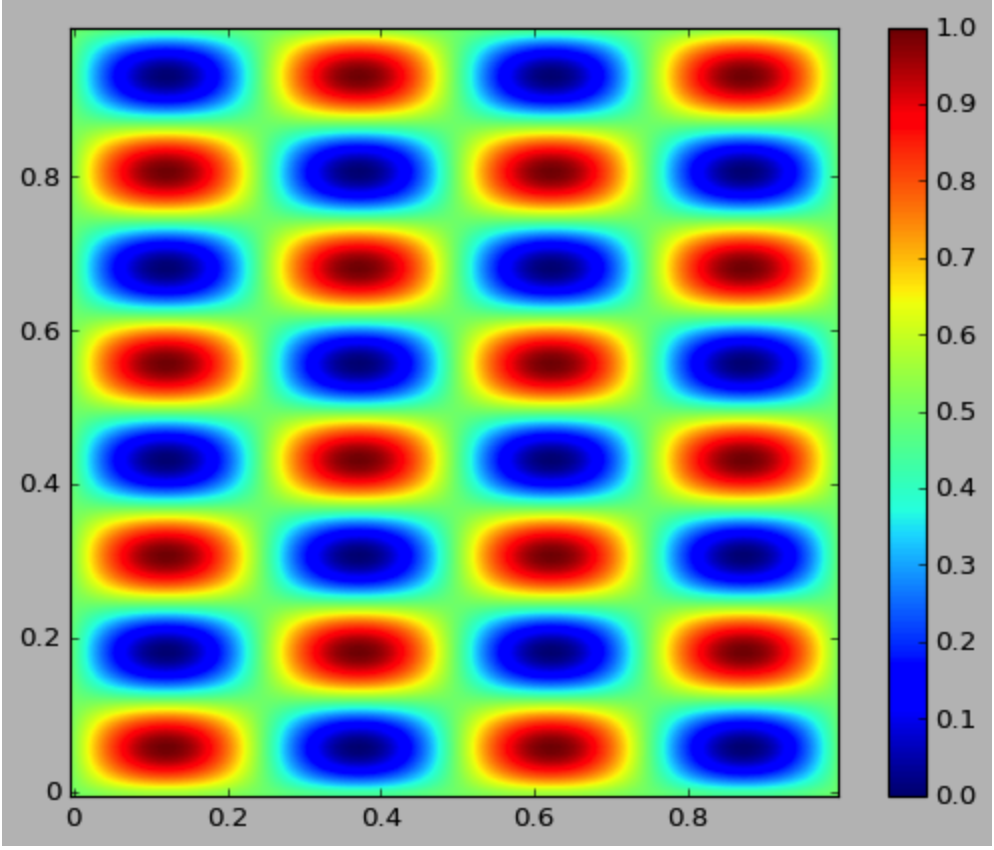}
\label{fig:f2}
\caption{Initial condition for Example 3.}
\end{figure}

In this example, we test a more complicated initial condition $S_0(x,y) = (1+\sin(4\pi x)\sin(8\pi y))/2$. Figure~\ref{fig:f2} shows the figure of  $S_0(x,y)$. This new initial condition is also bounded between $0$ and $1$. For this third numerical example, we will also consider the time step as $\Delta t=0.1$ and use the second example's velocity  given by
$$v(x,y)=(1+\sin^2(2\pi x)\sin(2\pi y),1+\sin(4\pi x))\cos(2\pi y))$$
 to test our algorithm. 
 The coarse and the fine meshes
are again $20\times 20$ and $100\times 100$ respectively. Same as before, we choose the condition $\|\bar{S}^{(m+1)}-\bar{S}^{(m)}\|\leq 0.0001$ as the stopping criterion.

\begin{figure}
\centering
  \includegraphics[scale=0.4]{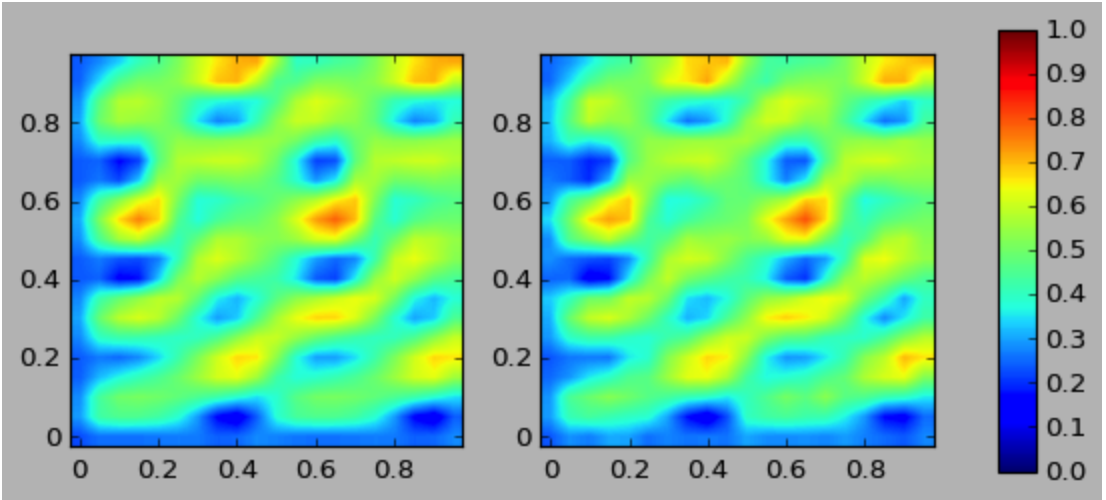}
    \includegraphics[scale=0.4]{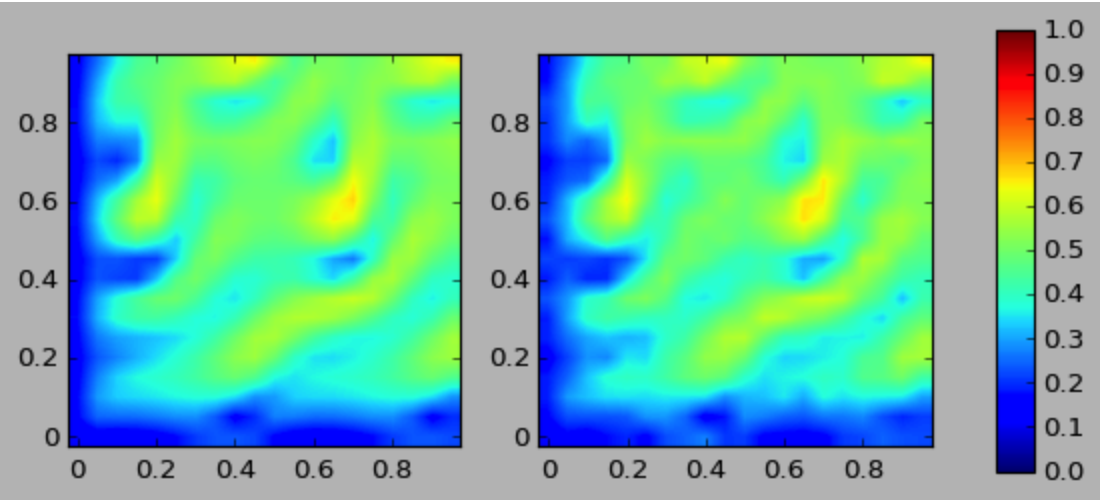}
      \includegraphics[scale=0.4]{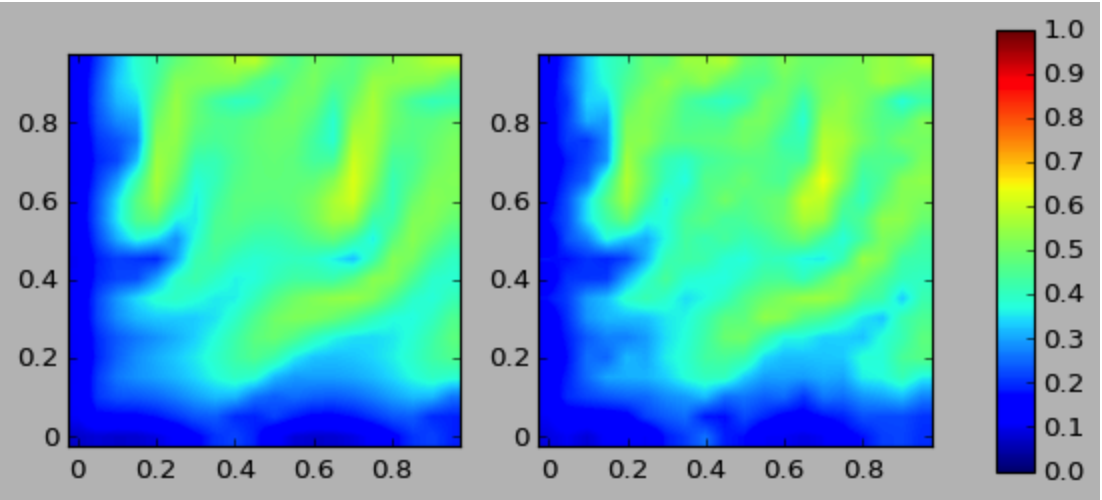}
        \includegraphics[scale=0.4]{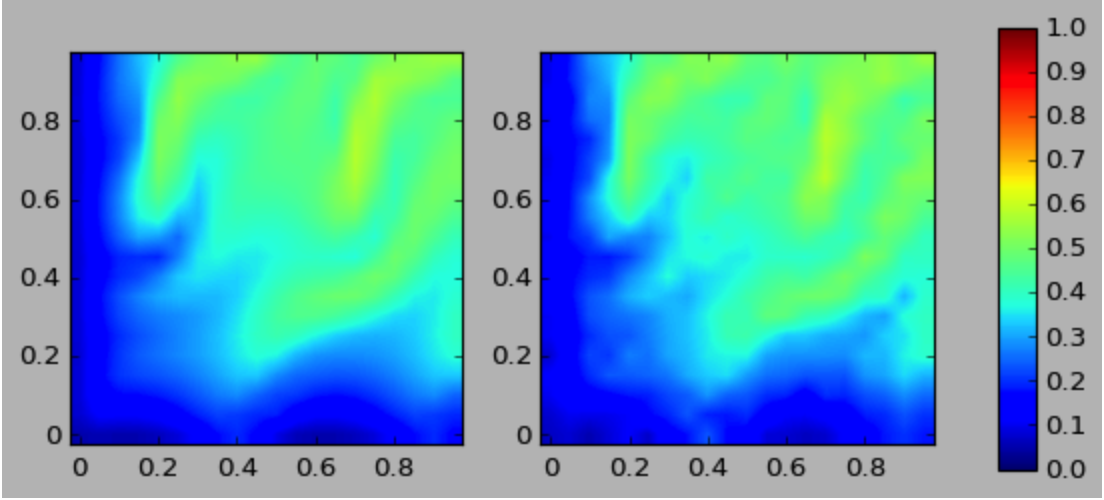}
          \includegraphics[scale=0.4]{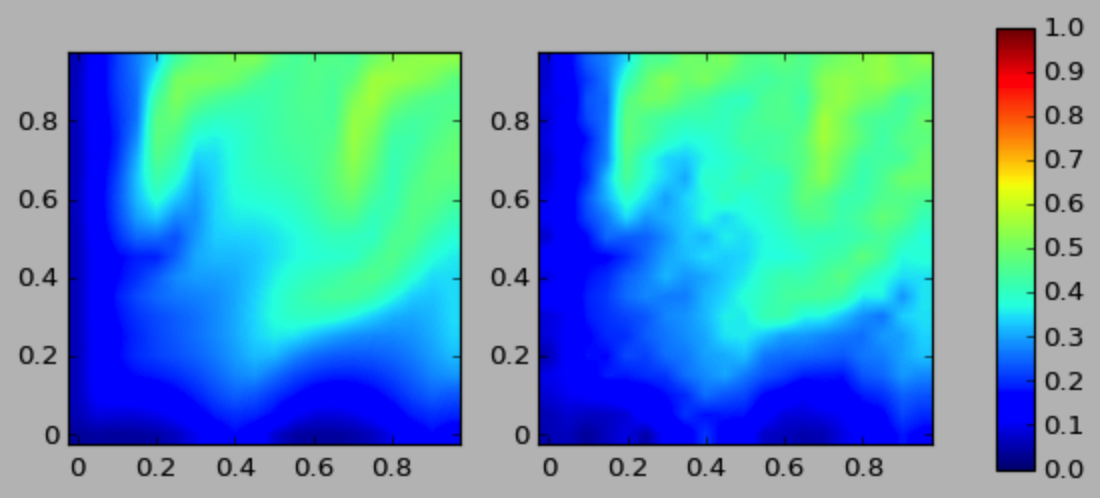}
\caption{Solutions at the first 5 time steps $t=0.1,0.2, \cdots, 0.5$ for Example 3. Left: upscaled reference solution. Right: upscaling numerical solution using deep learning based upscaling.}
\label{fig:solutionpic3}
\end{figure}

In this numerical example, we test the algorithm with 5 time steps.  Figure~\ref{fig:solutionpic3} show the first 5 time steps computed using our proposed scheme with deep neural network model. The corresponding average values of the reference solution on the coarse grid are shown in the same figure. We observe that these two solutions have very good agreement. 

\begin {table}[H]
 \caption{Root mean square error for Example 3 using neural network model}
  \label{fig:table3.1}
\begin{center}
 \begin{tabular}{|c|c|} 
 \hline
time step & relative error of the mean value \\ [0.5ex] 
 \hline\hline
$t=0.1$ & 0.03850  \\ 
 \hline
$t=0.2$&0.04167  \\ 
 \hline
$t=0.3$&0.04516  \\ 
 \hline
$t=0.4$&0.04901  \\ 
 \hline
 $t=0.5$&0.05219  \\ 
 \hline
\end{tabular}
\end{center}
\end {table}

For the proposed neural network model, the relative error of the numerical solution of first time step is 0.03850. The relative error of the numerical solution of fifth time step is 0.05219. Table~\ref{fig:table3.1} shows the relative errors of the numerical solution of using our proposed upscaling scheme with deep
neural network model. Also, we observe that when the time step increases, the relative error of the numerical solution rises slowly. Also the time used is about 50s for computation of one time step of using neural network model.  As a comparison, if we compute the local downscaling functions without the use of deep neural network, we need about 6000s for computation of one time step. Since the velocity of Examples 2 and 3 is the same, we only need to use the same neural network model from Example 2 and apply for this new initial condition.

\begin {table}[H]
 \caption{Root mean square error for Example 3 using the upscaled model without neural network}
  \label{fig:table3.2}
\begin{center}
 \begin{tabular}{|c|c| c|} 
 \hline
time step  & relative error of the mean value \\ [0.5ex] 
 \hline\hline
$t=0.1$ & 0.01859  \\ 
 \hline
$t=0.2$ & 0.02939  \\ 
 \hline
$t=0.3$ & 0.03762  \\ 
 \hline
$t=0.4$ & 0.04591  \\ 
 \hline
 $t=0.5$ & 0.05299  \\ 
 \hline
\end{tabular}
\end{center}
\end {table}

As a comparison, we report the results using the downscaling functions without neural network. The relative error of the numerical solution at the first time step is 0.01859. The relative error of the numerical solution at the fifth time step is 0.05499. Table~\ref{fig:table3.2} shows the relative errors of the numerical solution when our method is applied without the use of deep neural network. Also, the performance of solving the local downscaling functions numerically without the use of deep neural network and solving the local downscaling functions using neural network model are similar.

\begin {table}[H]
 \caption{Root mean square error for Example 3 using the finite volume method}
  \label{fig:table3.3}
\begin{center}
 \begin{tabular}{|c|c| c|} 
 \hline
time step  & relative error of the mean value \\ [0.5ex] 
 \hline\hline
$t=0.1$& 0.09739  \\ 
 \hline
$t=0.2$ & 0.11376  \\ 
 \hline
$t=0.3$ & 0.13648  \\ 
 \hline
$t=0.4$& 0.11288  \\ 
 \hline
 $t=0.5$ & 0.11446  \\ 
 \hline
\end{tabular}
\end{center}
\end {table}

The final comparison is the standard coarse-grid finite volume method with upwind flux. The relative error of the numerical solution at the first time step is 0.14304. The relative error of the numerical solution at fifth time step is 0.11446. Table~\ref{fig:table3.3} shows the relative errors of the numerical solution computed using the finite volume method.  We notice that our proposed nonlinear upscaling is able to give better numerical solutions.

\section{Conclusion}
\label{sec:conclusion}

In this paper, we develop deep learning based nonlinear upscaling for nonlinear transport equations with heterogeneous velocity.
The technique is based on the recently developed NLMC method. To construct the coarse scale model, local downscaling operations
are used to reconstruct fine scale information from coarse grid averages. Due to the nonlinearity of the problem, these local solutions are expensive to compute. 
To overcome this bottleneck, we propose the use of deep neural network to approximate this procedure. 
Since reduced model is used, the learning process is robust, and the resulting neural network gives accurate approximations. 
Our numerical results show promises of the proposed approach.

\section*{Acknowledgements}
Eric Chung's work is partially supported by the Hong Kong RGC General Research Fund (Project numbers 14304719 and 14302018) and the CUHK Faculty of Science Direct Grant 2018-19.
The first two authors thank NVIDIA - NVAITC for providing computing resources.

\bibliographystyle{plain}
\bibliography{references0,references11,references2}

\end{document}